\documentclass[reqno,a4paper,12pt]{amsart} 
    
    \usepackage{amsmath,amscd,amsfonts,amssymb}
    \usepackage{mathrsfs,dsfont}
    \usepackage{color}
    \usepackage{mathtools}
    \usepackage{hyperref}
    \usepackage{tikz-cd}

    \usepackage{tikz}
    \usetikzlibrary{decorations.pathreplacing}
    
    \numberwithin{equation}{section}
    \numberwithin{figure}{section}
    
    \def\R{\mathbb{R}}

    \def\Z{\mathbb{Z}}

    \def\S{\mathcal{S}}

    \def\Lam{\Lambda}

    \renewcommand\leq{\leqslant}
    \renewcommand\geq{\geqslant}

    \renewcommand{\ae}{{\mathrm{a.e.}}}

    \newcommand{\Mod}[1]{\ (\mathrm{mod}\ #1)}

    \theoremstyle{plain}
    \newtheorem{thm}{Theorem}[section]
    \newtheorem{theorem}[thm]{Theorem}
    
    \newtheorem{lemma}[thm]{Lemma}
    \newtheorem{corollary}[thm]{Corollary}

    \newtheorem{conjecture}[thm]{Conjecture}
    
    \newtheorem*{claim*}{Claim}
    
    \newcommand{\thmref}[1]{Theorem~\ref{#1}}

    \newcommand{\defref}[1]{Definition~\ref{#1}}

    \newcommand{\corref}[1]{Corollary~\ref{#1}}
    \newcommand{\conjref}[1]{Conjecture~\ref{#1}}
    \newcommand{\exampref}[1]{Example~\ref{#1}}

    \newcommand{\figref}[1]{Figure~\ref{#1}}
    
    \theoremstyle{definition}
    \newtheorem{definition}[thm]{Definition}
    \newtheorem*{definition*}{Definition}
    \newtheorem*{remarks*}{Remarks}
    \newtheorem*{remark*}{Remark}
    \newtheorem{remark}[thm]{Remark}
    \newtheorem{example}[thm]{Example}

    \newenvironment{enumerate-math}
    {\begin{enumerate}
    		\addtolength{\itemsep}{5pt}
    		}
    	{\end{enumerate}}

    \newenvironment{enumerate-text}
    {\begin{enumerate}
    		\addtolength{\itemsep}{5pt}
    		}
    	{\end{enumerate}}

\begin{document}

    	\title[A counterexample to the periodic tiling conjecture]{A counterexample to the periodic tiling conjecture (announcement)}
    
    	\author{Rachel Greenfeld}
    	\address{School of Mathematics, Institute for Advanced Study, Princeton, NJ 08540.}
    	\email{greenfeld.math@gmail.com}
    	\author{Terence Tao}
    	\address{UCLA Department of Mathematics, Los Angeles, CA 90095-1555.}
    	\email{tao@math.ucla.edu}

    	\subjclass[]{}
    	\date{}
    	
    	
    	\begin{abstract}  The periodic tiling conjecture asserts that any finite subset of a lattice $\Z^d$ which tiles that lattice by translations, in fact tiles periodically.  We announce here a disproof of this conjecture for sufficiently large $d$, which also implies a disproof of the corresponding conjecture for Euclidean spaces $\R^d$.  In fact, we also obtain a counterexample in a group of the form $\Z^2 \times G_0$ for some finite abelian $G_0$.  Our methods rely on encoding a certain class of ``$p$-adically structured functions'' in terms of certain functional equations.
    	\end{abstract}
    	
    	\maketitle
    	
    	\section{introduction}

        Let $G = (G,+)$ be a discrete abelian group.  If $A, F$ are subsets of $G$, we write $A \oplus F = G$ if every element of $G$ has a unique representation in the form $a+f$ with $a \in A$, $f \in F$. If this occurs, we say that $F$ \emph{tiles $G$ by translations}.  If in addition $A$ is \emph{periodic}, by which we mean that it is the finite union of cosets of a finite index subgroup of $G$, we say that $F$ \emph{periodically tiles $G$ by translations}.
        
        We can then pose the following conjecture:

        \begin{conjecture}[Discrete periodic tiling conjecture]\label{ptc}  Let $F$ be a finite non-empty subset of a finitely generated discrete abelian group $G$.  If $F$ tiles $G$ by translations, then $F$ periodically tiles $G$ by translations.
        \end{conjecture}

As observed by Wang \cite{wang}, this conjecture would imply as a corollary that the question of whether a given finite non-empty subset $F$ of some explicitly presented finitely generated abelian group (e.g., $\Z^d$) tiles that group by translations is (algorithmically) decidable (and logically decidable for any fixed $F$); see \cite[Appendix A]{GT2} for some further discussion. However, we will not discuss issues of decidability further here.

        We also phrase the following continuous analogue of this conjecture.  If $\Sigma$ is a bounded measurable subset of a Euclidean space $\R^d$, and $\Lam$ is a discrete subset of $\R^d$, we write $\Sigma \oplus \Lam =_\ae \R^d$ if the translates $\Sigma + \lambda$, $\lambda \in \R^d$, partition $\R^d$ up to null sets.  If this occurs, we say that $\Sigma$ (measurably) \emph{tiles $\R^d$ by translations}.  If in addition $\Lam$ is \emph{periodic}, by which we mean that it is the finite union of cosets of a lattice in $\R^d$, we say that $\Sigma$ (measurably) \emph{periodically tiles $\R^d$ by translations}.

\begin{conjecture}[Continuous periodic tiling conjecture]\label{ptc-cts}  Let $\Sigma$ be a bounded measurable subset of $\R^d$.  If $\Sigma$ tiles $\R^d$ by translations, then $\Sigma$ periodically tiles $\R^d$ by translations.
\end{conjecture}

The following partial results towards these conjectures are known:

\begin{itemize}
    \item Conjecture \ref{ptc} is trivial when $G$ is a finite abelian group, since in this case all subsets of $G$ are periodic.
\item  Conjectures \ref{ptc} and \ref{ptc-cts} were established in $G=\Z$ and $G=\R$ \cite{N,LW}.  The argument in \cite{N} also extends to the case $G = \Z \times G_0$ for any finite abelian group $G_0$ \cite[Section 2]{GT2}. 
\item When $G=\Z^2$, \conjref{ptc}  was established by Bhattacharya \cite{BH} using ergodic theory methods.  In  \cite{GT} we gave an alternative proof of this result, and  prove  tiling in $\Z^2$ is \textit{weakly periodic} (a disjoint union of finitely many one-periodic sets). 
\item When $G=\R^2$, \conjref{ptc-cts} is known to hold for any tile which is a topological disc  \cite{bn,gbn,ken}.
\item For $d>2$, the conjecture is known when the cardinality $|F|$ of $F$ is prime or equal to $4$ \cite{szegedy}, but remained open in general. 
       \item The continuous periodic tiling conjecture in $\R^d$ implies the discrete periodic tiling conjecture in $\Z^d$.  This implication is standard, arising from ``encoding'' a discrete subset $F$ of $\Z^d$ as a bounded measurable subset $F \oplus \Omega$ in $\R^d$, where $\Omega$ is a ``generic'' fundamental domain of $\R^d/\Z^d$; we provide the details in our forthcoming paper \cite{greenfeld-tao-ptc}.
    \item In \cite{bgu}, it was recently shown that the discrete periodic tiling conjecture in $\Z^d$ also implies the discrete periodic tiling conjecture in every quotient group $\Z^d/\Lambda$.
    \item The analogues of the above conjectures are known to fail when one has two or more translational tiles instead of just one; see \cite{GT2} (particularly Table 1) for a summary of results in this direction.  In particular, in \cite[Theorems 1.8, 1.9]{GT2} it was shown that the analogue of  \conjref{ptc} for two tiles fails\footnote{Strictly speaking, the counterexample in that paper involved tiling a periodic subset $E$ of the group $G$, rather than the full group $G$.} for $\Z^2 \times G_0$ for some finite group $G_0$, and also for $\Z^d$ for some $d$.
\end{itemize}

In our forthcoming paper \cite{greenfeld-tao-ptc} we will obtain counterexamples to the above conjecture.  Our main result is

\begin{theorem}[Counterexample to \conjref{ptc}]\label{main}  There exists a finite group $G_0$ such that the discrete periodic tiling conjecture fails for $\Z^2 \times G_0$.
\end{theorem}

The group $\Z^2 \times G_0$ can be viewed as a quotient $\Z^d/\Lambda$ of a lattice $\Z^d$, so by the preceding implications, we obtain

\begin{corollary}[Counterexample to \conjref{ptc-cts}]\label{main-cor}  For sufficiently large $d$, the discrete periodic tiling conjecture fails for $\Z^d$, and the continuous periodic tiling conjecture fails for $\R^d$.
\end{corollary}

Our methods produce a finite group $G_0$, and hence a dimension $d$, that is in principle explicitly computable, but we have not attempted to optimize the size of these objects.  In particular the dimension $d$ produced by our construction will be extremely large.

On the other hand, \conjref{ptc} is known to be false for multiple tiles $F_1,\dots,F_k$ in $\Z^2$. In previous literature,  aperiodic \emph{multi-tilings}  have been created by such devices as the ``cut and project'' method of Meyer \cite{M70,M95,L}, the finite state machine approaches of Kari and Culik \cite{kari,culik}, or by encoding arbitrary Turing machines into a tiling problem \cite{Ber,Ber-thesis,ollinger}, and other methods (see \cite{S} for a survey on aperiodic tilings and constructions).   Unfortunately, we were not able to adapt any of these methods to the setting of a \emph{single} translational tile.  Instead, our source of aperiodicity is more novel, in that our tiling of $\Z^2 \times G_0$ is forced to exhibit a ``$p$-adic'' structure  for some large enough but fixed prime $p$, say $p>48$, in the sense that for each power $p^j$ of $p$, the tiling is periodic with period $p^j \Z^2 \times \{0\}$ outside of a small number of cosets of that subgroup $p^j \Z^2 \times \{0\}$, but is unable to be genuinely periodic with respect to any of these periods.  To achieve this we will set up a certain ``Sudoku-type puzzle'', which will be rigid enough to force all solutions of this problem to exhibit $p$-adic (and therefore non-periodic) behavior, yet is not so rigid that there are no solutions whatsoever.  By modifying arguments from our previous paper \cite{GT2}, we are then able to encode this Sudoku-type puzzle as an instance of the original tiling problem $A \oplus F = \Z^2 \times G_0$.

In the remainder of this announcement, we sketch the main steps involved in establishing  \thmref{main}.  Full details of these steps will appear in the forthcoming paper \cite{greenfeld-tao-ptc}; we give a brief summary here, and then expand upon the key steps in the remaining sections of this announcement.

Our argument is a variant of the construction used in our previous paper \cite{GT2} to produce counterexamples to the periodic tiling conjecture for two tiles, although the fact that we are now tiling the whole group $G$ instead of a periodic subset of $G$, and that we are only allowed to use one tile instead of two, creates additional technical challenges; see, e.g., the discussion in our previous paper  \cite[Section 10]{GT2}.

As in \cite{GT2}, we begin by replacing the single tiling equation $A \oplus F = G$ with a system $A \oplus F^{(m)} = G$, $m = 1,\dots,M$ of tiling equations for an arbitrary $M$, by an elementary ``stacking'' procedure that takes advantage of our freedom to enlarge the group $G$.  This creates a flexible ``tiling language'' of constraints on the tiling set $A$; the challenge is to use this language to obtain a system of constraints that is strict enough to force aperiodic behavior on this set $A$, while simultaneously being relaxed enough to admit at least one solution.

Next, we again follow \cite{GT2} and pass from this tiling language to a more familiar language of functional equations, basically by spending one of the equations $A \oplus F^{(m)} = G$ in the system to force the tiling set $A$ to be a graph of a function $f=(f_1,\dots,f_K)$, where  $f_i \colon \Z^2 \times G_0 \to \Z/2q\Z$, $1\leq i\leq K$, and  $G_0$ is an additional small finite abelian group which we retain for technical reasons.  One can then use one or more tiling equations $A \oplus F^{(m)} = G$ in the tiling language to enforce useful functional constraints on these functions $f_i$.  For instance, one can ensure that a given function $f_i$ exhibits periodicity in some direction $v_i \in \Z^2$, or that it encodes the $\mod N$ reduction map $n \mapsto n \mod N$ (up to a shift) for a given modulus $N$.  Crucially, we are also able to encode the assertion that a certain subcollection of the $f_i$ (after a routine normalization) take values in the two-element set $\{-1,1\} \mod 2q\Z$, thus effectively making them boolean functions.  By adapting a construction from \cite[Section 7]{GT2}, we can then use tiling equations to encode arbitrary pointwise constraints
\begin{equation}\label{fkx}
(f_1(x), \dots, f_K(x)) \in \Omega
\end{equation}
for all $x \in \Z^2 \times G_0$ and arbitrary subsets $\Omega$ of $\{-1,1\}^K$.  

By some further elementary transformations, we are then able to reduce matters to demonstrating aperiodicity of a certain ``Sudoku-type puzzle''.  In this puzzle, we have an unknown function $f \colon \{1,\dots,N\} \times \Z \to (\Z/p\Z)^\times$ on a vertically infinite ``Sudoku board'' $\{1,\dots,N\} \times \Z$ which fills each cell $(n,m)$ of this board with an element $f(n,m)$ of the multiplicative group $(\Z/p\Z)^\times$ for some fixed but large prime $p$ (for instance, one can take $p=53$).  Along every row or diagonal (and more generally along any non-vertical line) of this board, the function $f$ is required\footnote{This is analogous to how, in the most popular form of a Sudoku puzzle, the rows, columns, and $3 \times 3$ blocks of cells on a board $\{1,\dots,9\} \times \{1,\dots,9\}$ are required to be permutations of the digit set $\{1,\dots,9\}$.} to exhibit ``$p$-adic behavior''; the precise description of this behavior will be given below, but roughly speaking we will require that on each such non-vertical line, $f$ behaves like a rescaled version of the function \begin{equation}\label{fpn}
f_p(n) \coloneqq \frac{n}{p^{\nu_p(n)}} \mod p
\end{equation}
(where $\nu_p(n)$ is the number of times $p$ divides $n$), 
that assigns to each integer $n$ the final non-zero digit in its base $p$ expansion (with the convention $f_p(0) \coloneqq 1$).  We also impose a non-degeneracy condition that the function $f$ is not constant along any of its columns.  For suitable choices of parameters $p,N$, we will be able to ``solve'' this Sudoku problem and show that solutions to this problem exist, but necessarily exhibit $p$-adic behavior along the infinite columns of this puzzle, and in particular are non-periodic.  By combining this aperiodicity result with the previous encodings and reductions, we are able to establish  \thmref{main} and hence \corref{main-cor}.

    	\subsection{Acknowledgments}
    	RG was partially supported by NSF grant DMS-2242871. TT was partially supported by NSF grant DMS-1764034 and by a Simons Investigator Award.

\section{Equivalence with the multiple periodic tiling conjecture}\label{sec:m-ptc}

\conjref{ptc} involves a single tiling equation $A \oplus F = G$.  It turns out that the following conjecture, despite seeming stronger than \conjref{ptc} as it involves multiple tiling equations $A \oplus F^{(m)} = G$), in fact, is equivalent  to \conjref{ptc}:

\begin{conjecture}[Multiple periodic tiling conjecture]\label{ptc-multi}  Let $F_1,\dots,F_M$ be finite non-empty subsets of a finitely generated discrete abelian group $G$.  If there exists a solution $A \subset G$ to the system of tiling equations $A \oplus F^{(m)} = G$ for $m=1,\dots,M$, then there exists a \emph{periodic} solution $A' \subset G$ to the same system $A' \oplus F^{(m)} = G$ for $m=1,\dots,M$.
\end{conjecture}

\begin{figure}
    \centering
    \includegraphics[width = .4\textwidth]{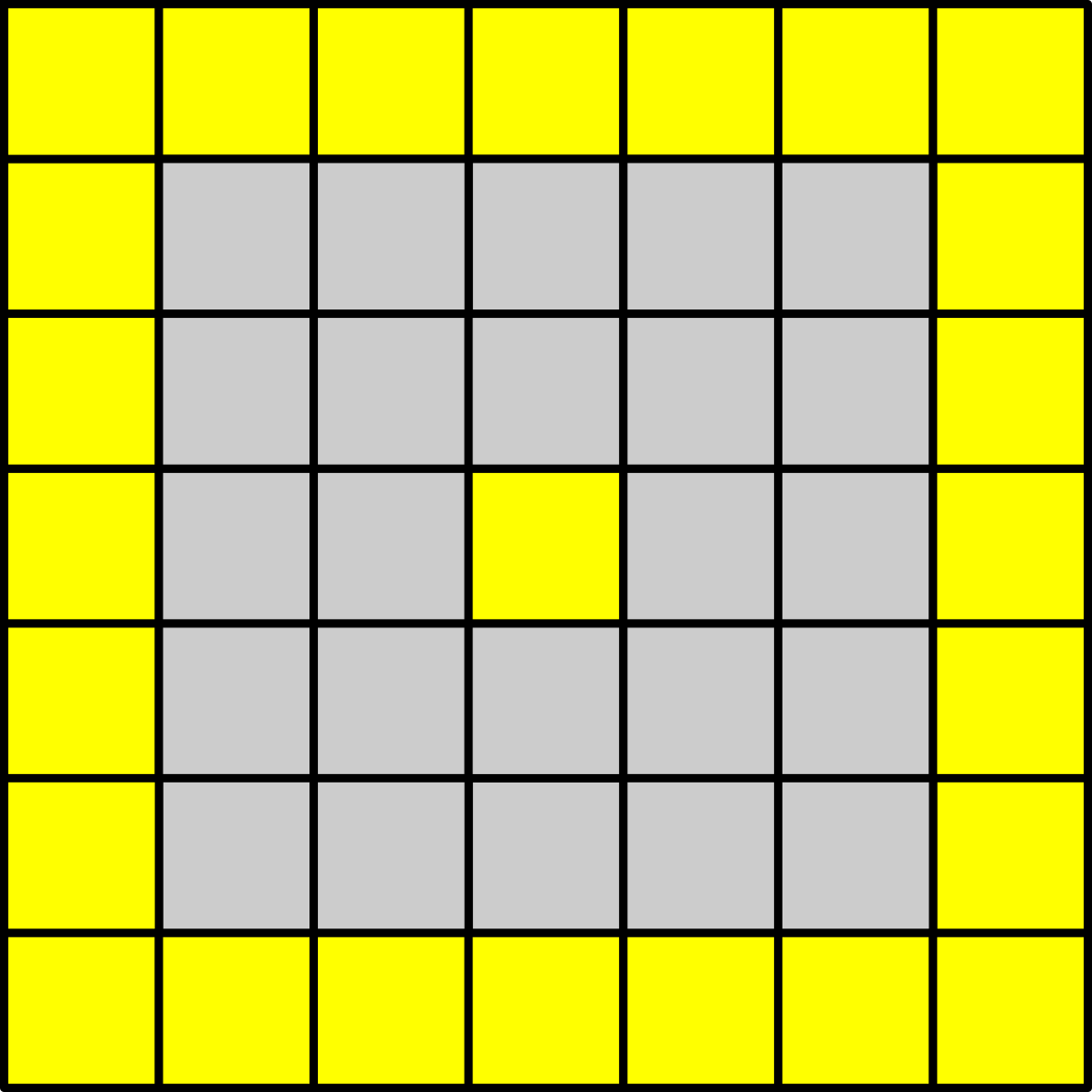}
    \caption{For $M=2$, we partition  $H=(\Z/7\Z)^2$ into $E_1$ (yellow) and $E_2$ (grey), such that \eqref{intersective} is satisfied.} 
    \label{fig:partition}
\end{figure}

Clearly \conjref{ptc} is the special case $M=1$ of \conjref{ptc-multi}.  In the converse direction, one can use the freedom to enlarge the group $G$ to deduce \conjref{ptc-multi} for general $M$ from   \conjref{ptc}.  This sort of implication first appeared in \cite[Theorem 1.15]{GT2} in a slightly different setting, in which one only sought to tile a periodic subset $E$ of the group $G$.  However, it turns out that the arguments can be adapted to the setting in which one is tiling the entire group $G$.  To illustrate the ideas let us just focus on the $M=2$ case.  More precisely, we are given finite non-empty subsets $F^{(1)}$, $F^{(2)}$ of $G$ for which there is a solution $A \subset G$ to the system
\begin{equation}\label{afg}
 A \oplus F^{(1)} = A \oplus F^{(2)} = G,
 \end{equation}
and we wish to somehow use   \conjref{ptc} to then construct a periodic solution $A' \subset G$ to the same system
\begin{equation}\label{apfg}
A' \oplus F^{(1)} = A' \oplus F^{(2)} = G.
\end{equation}
It turns out that in the finite abelian group $H \coloneqq (\Z/7\Z)^2$ one can find a partition $H = E_1 \uplus E_2$ in such a way that one has
\begin{equation}\label{intersective} (E_i + h_i) \cap (E_j + h_j) \neq \emptyset
\end{equation}
for any $i,j \in \{1,2\}$ and $h_i,h_j \in H$, unless $i \neq j$ and $h_i=h_j$; see   \figref{fig:partition} for an explicit example.  Informally, this means that the only way to partition $H$ in terms of translates of $E_1, E_2$ is by the partitions $H = (E_1 + h) \uplus (E_2 + h)$ for $h \in H$.

Now we work in the product group $G \times H$ and introduce the single tile
$$ \tilde F \coloneqq (F^{(1)} \times E_1) \uplus (F^{(2)} \times E_2).$$
It is not hard to see that the system \eqref{afg} implies the single tiling equation
$$ (A \times \{0\}) \oplus \tilde F = G \times H.$$
Applying   \conjref{ptc} to the product group $G \times H$, we conclude that there exists a periodic subset $\tilde A'$ of $G \times H$ such that
\begin{equation}\label{atile}
 \tilde A' \oplus \tilde F = G \times H.
 \end{equation}
One can use the intersective property \eqref{intersective} to conclude that $\tilde A'$ intersects each vertical fiber $\{g\} \times H$ of $G \times H$, $g\in G$, in at most one point, and hence $\tilde A'$ must be a graph
$$ \tilde A' = \{ (a, f(a)): a \in A' \}$$
for some periodic $A' \subset G$ and some function $f \colon A' \to H$.  One can then use the tiling equation \eqref{atile} (and the intersective property \eqref{intersective}) to conclude the pair of tiling equations \eqref{apfg}, thus establishing the $M=2$ case of  \conjref{ptc-multi}.  The case of general $M$ can be established by a similar argument; see \cite{greenfeld-tao-ptc} for details.

\section{Reduction to constructing an aperiodic system of functional equations}\label{sec:sfe}

Our task is now to construct a system of tiling equations $A \oplus F^{(m)} = G$ that admits solutions, but no periodic solutions.  It will be easier to work with functional equations rather than tiling equations.  For this we make a basic observation: if $A$ is a subset of a product group $G \times H$, where $H$ is some finite abelian group, then $A$ obeys the tiling equation
$$ A \oplus (\{0\} \times H) = G \times H$$
if and only if $A$ is the graph 
\begin{equation}\label{graph}
A = \{ (a, f(a)): a \in G \}
\end{equation}
of some function $f \colon G \to H$.  This gives us a correspondence between certain tiling sets $A \subset G \times H$ and functions $f \colon G \to H$ that allows us to translate tiling equations (or systems of tiling equations) as functional equations involving $f$.

For our construction, the range group $H$ will be taken to be of the form $(\Z/2q\Z)^K$ for some large prime $q$ and some natural number $K$.  The function $f \colon G \to H$ can then be thought of as $K$ independent functions $f_1,\dots,f_K \colon G \to \Z/2q\Z$ taking values in the cyclic group $\Z/2q\Z$.  We will then use a variety of functional equations (which can then be translated into tiling equations), to create suitable constraints between these functions, which will ultimately allow us to encode a certain ``Sudoku puzzle''.

The full translation procedure is rather intricate and requires a certain amount of notation, but we illustrate some aspects of this translation by means of simple examples in which constraints on a function $f \colon G \to H$ (or on multiple functions $f_i \colon G \to \Z/2q\Z$) are encoded as one or more tiling equations $A \oplus F^{(m)} = G \times H$, where $A$ is a graph \eqref{graph}.

\begin{example}[Encoding periodicity] Let $A \subset G \times H$ be a graph \eqref{graph} of a function $f \colon G \to H$, and let $v \in G$.  Then the periodicity property
$$ f(x+v) = f(x)\quad  \forall x \in G$$
is equivalent to the tiling equation
$$ A \oplus ((\{0\} \times \{0\}) \uplus (\{v\} \times (H \backslash \{0\}))) = G \times H.$$
\end{example}

\begin{example}[Encoding a shifted mod $N$ function]\label{shift}  Let $A \subset \Z \times \Z/N\Z$ be a graph \eqref{graph} of a function $f \colon \Z \to \Z/N\Z$.  Then the property that 
$$ f(x) = x + c \mod N$$
for all $x \in \Z$ and some (unspecified) shift $c \in \Z$ is equivalent to the tiling equation
$$ A \oplus ((\{0\} \times \{0\}) \uplus (\{1\} \times (\Z/N\Z \backslash \{1\}))) = \Z \times \Z/N\Z.$$
\end{example}

\begin{remark} The presence of the shift $c$ in the above example (and in several of the examples below) reflects a certain translation invariance in this problem: if the graph \eqref{graph} of a function $f \colon G \to H$ obeys a system of tiling equations $A \oplus F^{(m)} = G$, then the same is true for any shift $f+c$ of that function.  This translation invariance was not present in our previous work \cite{GT2}, in which we tiled a periodic subset of $G$ rather than all of $G$, and thus causes some technical difficulties in our arguments in \cite{greenfeld-tao-ptc}; however, in many situations we will be able to normalize such shifts, for instance to equal zero.
\end{remark}

\begin{example}[Encoding linear constraints]\label{linear}  Let $f_1,\dots,f_K \colon \Z \to \Z/2q\Z$ be functions, and let $A \subset \Z \times (\Z/2q\Z)^K$ be the graph \eqref{graph} of the combined function $f = (f_1,\dots,f_K) \colon \Z \to (\Z/2q\Z)^K$.  Let $a_1,\dots,a_K \in \Z/2q\Z$ be coefficients.  Then the property that a linear relation
$$ a_1 f_1(x) + \dots + a_K f_K(x) = c$$
holds for all $x \in \Z$ and some (unspecified) constant $c \in \Z/2q\Z$ is equivalent to the tiling equation
$$ A \oplus ((\{0\} \times E) \uplus (\{1\} \times (H \backslash E))) = \Z \times (\Z/2q\Z)^K,$$
where $E \leq (\Z/2q\Z)^K$ is the subgroup
$$ E \coloneqq \{ (y_1,\dots,y_K) \in (\Z/2q\Z)^K: a_1 y_1 + \dots + a_K y_K = 0 \}.$$
\end{example}

\begin{example}[Encoding a rescaled boolean function]\label{bool}  Let $A \subset \Z \times \Z/2\Z \times \Z/2q\Z$ be the graph \eqref{graph} of a function $f \colon \Z \times \Z/2\Z \to \Z/2q\Z$.  Then the property that there exists an even $a \in \Z/2q\Z$ and odd $b \in \Z/2q\Z$ such that
$$ \{ f(x,0), f(x,1) \} = \{a,b\}$$
for all $x \in \Z$ (i.e., one has the anti-symmetry $f(x,y+1) = a+b-f(x,y)$ for all $(x,y) \in \Z \times \Z/2\Z$, and $f$ only takes on the two values $a,b$), is equivalent to the system of tiling equations
$$ A \oplus F = \Z \times \Z/2\Z \times \Z/2q\Z$$
and
$$ A \oplus (F \backslash \{(0,0,0),(0,1,0)\} \cup \{(1,0,0),(1,1,0)\}) = \Z \times \Z/2\Z \times \Z/2q\Z$$
where 
$$ F \coloneqq \{ (0,0), (0,1) \} \times 2\Z/2q\Z.$$
\end{example}

In the above example, one can use a translation invariance to normalize $a=0$, and then a dilation invariance (cf., the ``dilation lemma'' for translational tilings, established for instance in \cite[Lemma 3.1]{GT}) can also be used to normalize $b=1$.  However, when working with multiple functions $f_1,\dots,f_K \colon \Z \times \Z/2\Z \to \Z/2q\Z$, the analogous construction to  \exampref{bool} will end up placing the values of each $f_i$ in a different doubleton set $\{a_i,b_i\}$ for some even $a_i \in \Z/2q\Z$ and odd $b_i \in \Z/2q\Z$.  As before, the translation invariance will allow one to normalize all the $a_i$ to equal zero, but it turns out that the dilation invariance only permits us to normalize \emph{one} of the $b_i$ to equal $1$.  We will need to ``align'' these rescaled boolean functions so that all the spacings $|a_i-b_i|$ are equal, which will effectively allow us to normalize $a_i=0, b_i=1$ for all $i$.  To this end, the following elementary lemma turns out to be useful:

\begin{lemma}[Alignment lemma]  Let $G$ be a finitely generated abelian group.  For each $i=1,2,3,4$, let $f_i \colon G \to \{a_i,b_i \} \subset \Z/2q\Z$ be a function taking two values $a_i,b_i$ with $a_i-b_i$ odd.  Suppose that the function
$$ f_1 + f_2 - (2f_3 + f_4)$$
is constant on $G$, and the map $(f_1,f_2) \colon G \to \{a_1,b_1\} \times \{a_2,b_2\}$ is surjective.  Then we have $|a_1-b_1|=|a_2-b_2|$.
\end{lemma}

To illustrate this lemma, suppose that $f_1,f_2$ both take values in $\{0,1\}$; then $f_1+f_2$ takes values in $\{0,1,2\} \subset 2 \cdot \{0,1\} + \{0,1\}$, and so can be expressed as $f_1+f_2=2f_3+f_4$ for some functions $f_3,f_4 \colon G \to \{0,1\}$.  However, if $f_1$ takes values in $\{0,1\}$ and $f_2$ takes values in $\{0,3\}$, then $f_1+f_2$ now takes values in $\{0,1,3,4\}$ (with all four values attained if we assume $(f_1,f_2)$ surjective) and there is no obvious way to express this function as $2f_3+f_4$ where $f_3,f_4$ each take on just two values differing by an odd number.

After an appropriate use of this lemma and some renormalizations, we will be able to encode systems of linear equations involving boolean functions $$f_1,\dots,f_K \colon G \to \{-1,+1\}.$$  As observed in our previous paper \cite[Sections 6 and 7]{GT2}, pointwise constraints on such functions can be encoded using linear equations between these functions, which can in turn be encoded as tiling equations using variants of  \exampref{linear}.  We illustrate this with a simple example:

\begin{example}[Encoding a boolean constraint]  Let $f_1,f_2,f_3 \colon G \to \{-1,+1\}$ be boolean functions.  Then one has the relation
$$ (f_1(x), f_2(x), f_3(x)) \neq \{(-1,-1,-1),(+1,+1,+1)\}$$
for all $x \in G$ if and only if there exists another boolean function $f_4 \colon G \to \{-1,+1\}$ for which one has the linear relation
$$ f_1(x)+f_2(x)+f_3(x) = f_4(x).$$
\end{example}

With more complicated variants of this example one can encode arbitrary constraints of the form \eqref{fkx}; see \cite[Section 7]{GT2}.

\section{Reduction to demonstrating the aperiodicity of a ``Sudoku puzzle''}

\begin{figure}
    \centering
    \includegraphics[width = 1.1\textwidth]{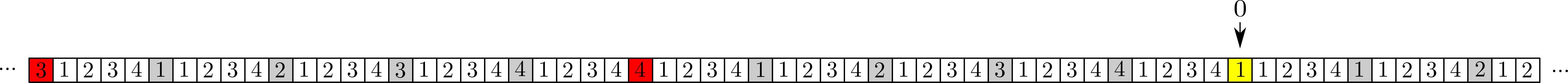}
    \caption{The function $f_p$ for $p=5$. The white cells correspond to $\nu_p(n)=0$, the grey cells are those with $\nu_p(n)=1$, the red ones have $\nu_p(n)=2$ and  yellow indicates $\nu_p(n)=+\infty$.} 
    \label{fig:fp}
\end{figure}

We now describe a certain ``Sudoku puzzle'' that can be encoded using the methods outlined in  sections \ref{sec:m-ptc} and \ref{sec:sfe} as a system of tiling equations.  We fix a large prime $p$ (say, larger than $48$) and set $N \coloneqq p^2$.  We introduce the ``$p$-adically structured function''
$f_p \colon \Z \to (\Z/p\Z)^\times$ by the formula \eqref{fpn} for $n \neq 0$, with $f_p(0) \coloneqq 1$. As mentioned in the introduction, $f_p(n)$ is simply the final non-zero digit of the base $p$ expansion of $n$, with the convention that the final non-zero digit of $0$ is $1$; see  \figref{fig:fp}.  Note that this function is not quite periodic; for each $j \geq 0$ it is a $p^j$-periodic function outside of a single coset $0 + p^j\Z$.

\begin{remark} The function $f_p$ can also be defined as the unique solution to the functional equations
\begin{align*}
    f_p(0) &= 1 \\
    f_p(n) &= n \mod p \hbox{ whenever } n \neq 0 \mod p \\
    f_p(pn) &= f_p(n) \hbox{ for all } n \in \Z.
\end{align*}
These are not quite the type of equations that can be encoded as tiling equations using the techniques sketched in the previous section - for instance, they fail to be translation-invariant - but, as these are quite basic functional equations, this does indicate some hope that one could try to capture the aperiodic ``$p$-adically structured'' nature of $f_p$, up to natural symmetries such as translation and dilation, using a suitable system of tiling equations.
\end{remark}

We now define a related class of $p$-adically structured functions that are localized rescaled versions of the base function $f_p$.

\begin{figure}
    \centering
    \includegraphics[width = .9\textwidth]{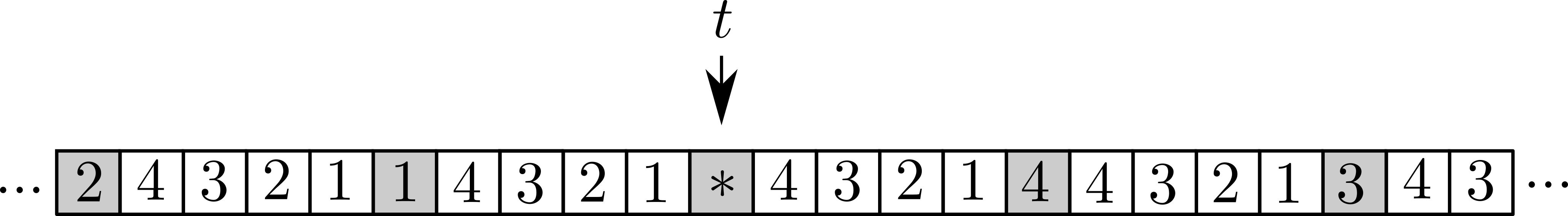}
    \caption{A non-constant element $g$ of $\S^2_p(\Z)$ with $p=5$ and step $h = 4 = -1 \Mod p$, depicted pictorially as an infinite string of cells, each filled with an element of $(\Z/p\Z)^\times$, which we identify with $\{1,2,3,4\}$ by a slight abuse of notation.  The function is affine outside of a ``bad coset'' $t + p\Z$, shaded in gray.  The asterisk depicts the fact that there is no constraint on $g$ in $t + p^2 \Z$.}
    \label{fig:badcoset}
\end{figure}

\begin{figure}
    \centering
    \includegraphics[width = .9\textwidth]{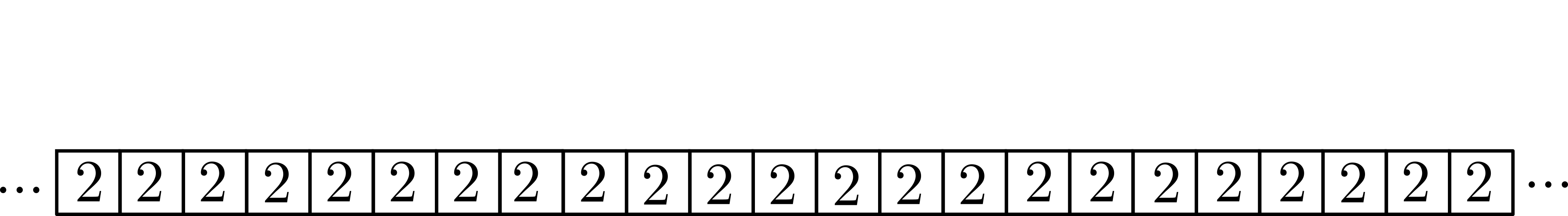}
    \caption{A constant element of $\S^2_p(\Z)$.  In this case, there is no ``bad coset''.}
    \label{fig:nobadcoset}
\end{figure}

\begin{definition}\label{s2p}  Let $\S^2_p(\{1,\dots,N\})$ denote the class of functions $g \colon \{1,\dots,N \} \to (\Z/p\Z)^\times$ which are either constant, or such that there exists $t \in \Z$ and $h \in (\Z/p\Z)^\times$ such that
\begin{equation}\label{key}
g(n)= h f_p(n-t)
\end{equation}
for all $n \in \{1,\dots,N\}$ with $n \neq t \mod p^2$; see Figures \ref{fig:badcoset}, \ref{fig:nobadcoset}.
\end{definition}

\begin{figure}
    \centering
    \includegraphics[width = .9\textwidth]{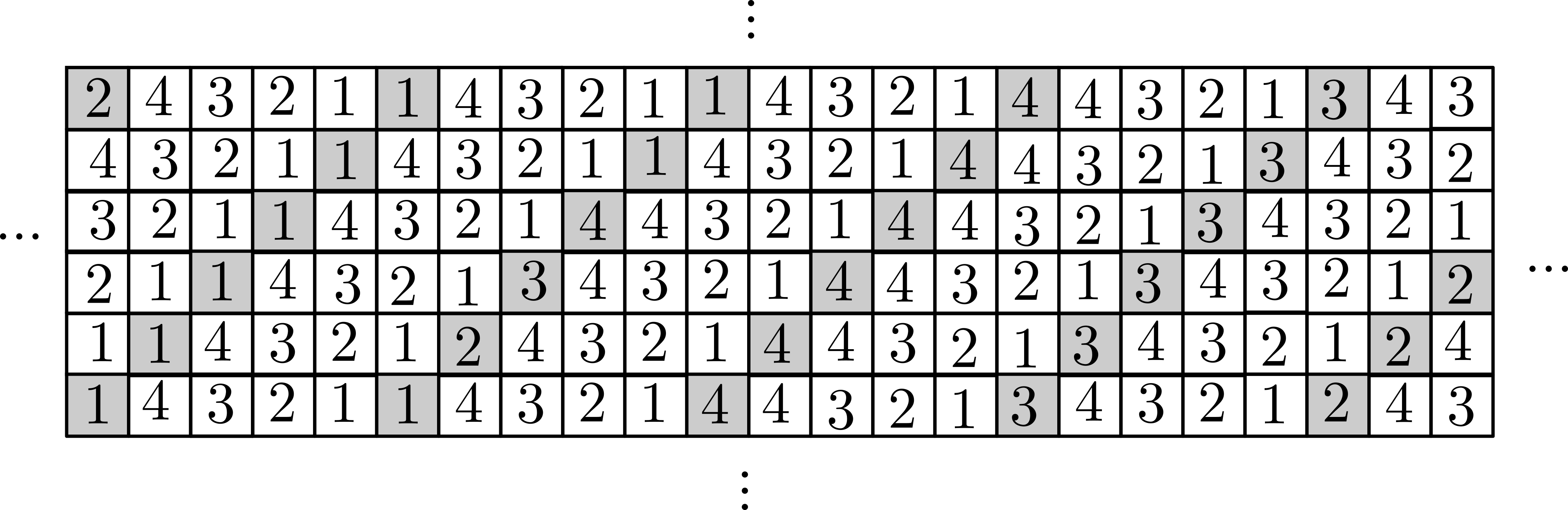}
    \caption{A portion of a $\S^2_5$-Sudoku solution.  Observe that it is affine outside of the shaded cells.}
    \label{fig:sudoku}
\end{figure}

Using this class of functions, we set up a ``Sudoku puzzle'':

\begin{definition}[$\S^2_p$-Sudoku]\label{Ssudoku}  A \emph{$\S^2_p$-Sudoku solution} is a function $F \colon \{1,\dots,N\} \times \Z \to (\Z/p\Z)^\times$ with the property that for every slope $j \in \Z$ and intercept $i \in \Z$, the function $n \mapsto F(n, jn+i)$ lies in the class $\S^r_p(\{1,\dots,N\})$. (See  \figref{fig:sudoku}.)
\end{definition}

A trivial example of a $\S^2_p$-Sudoku solution is a constant solution $F(n,m)=c$ for any constant $c \in (\Z/p\Z)^\times$.  A somewhat less obvious example is a function of the form $F(n,m) = f_p(an+bm+c)$ for some integers $a,b,c$: see \figref{fig:sudoku}.  Note that, due to the non-periodic nature of $f_p$, such solutions non-periodic, unless $b=0$ (in which case, the solution is constant on columns).  In \cite{greenfeld-tao-ptc} it is also convenient to consider the more general notion of an $\S^1_p$-Sudoku solution, in which the condition \eqref{key} is only enforced for $n \neq to \mod p$, but for simplicity of exposition we do not discuss this variant concept here.

A key component of the arguments in \cite{greenfeld-tao-ptc} is the proof of the following assertion:

\begin{theorem}[Aperiodicity of $\S^2_p$-Sudoku]\label{sudoku-aperiodic} Let $F \colon \{1,\dots,N\} \times \Z \to (\Z/p\Z)^\times$ be a $\S^2_p$-Sudoku solution, such that no column of $F$ is constant.  Then no column of $F$ is periodic.
\end{theorem}

We will discuss the proof of   \thmref{sudoku-aperiodic} in the next section.  For now, we indicate how this result, when combined with the encoding techniques from the previous section, permit a counterexample to  \conjref{ptc-multi} and hence Conjectures \ref{ptc}, \ref{ptc-cts}.

Given a $\S^2_p$-Sudoku solution $F \colon \{1,\dots,N\} \times \Z \to (\Z/p\Z)^\times$, one can define the functions $f_n \colon \Z^2 \to (\Z/p\Z)^\times$ for $n=1,\dots,N$ by the formula
\begin{equation}\label{fnij}
f_n(i,j) \coloneqq F(n, jn+i)
\end{equation}
for all $i,j \in \Z$.  From the definition of a $\S^2_p$-Sudoku solution, one can verify the following properties:

\begin{itemize}
    \item[(i)]  Each function $f_n$ obeys the periodicity property $f_n(i,j) = f_n(i-n,j+1)$ for all $(i,j) \in \Z^2$.
    \item[(ii)]  For each $(i,j)$, the tuple $(f_1(i,j),\dots,f_N(i,j))$ lies in $\S^2_p$.
\end{itemize}

Conversely, given a collection of functions $f_1,\dots,f_N \colon \Z^2 \to (\Z/p\Z)^\times$ obeying the axioms (i), (ii), one can find a unique $\S^2_p$-Sudoku solution  $F \colon \{1,\dots,N\} \times \Z \to (\Z/p\Z)^\times$ obeying the property \eqref{fnij}.  Thus, the property of being a $\S^2_p$-Sudoku can be encoded by the type of constraints discussed in the previous section; and after some further effort can therefore be encoded as a system of tiling equations; see \cite{greenfeld-tao-ptc} for details.

If one takes the non-periodic $\S^2_p$-Sudoku solution $F(n,m) = f_p(m)$, we see for instance that $f_1(i,j) = i+j \mod p$ whenever $i+j \neq 0 \mod p$.  Using a variant of   \exampref{shift}, this type of constraint (up to a translation) can also be encoded as a tiling equation.  Because of this, we can generate a system of tiling equations that generates solutions, all of whom are necessarily non-periodic, which will lead to our proof of  \thmref{main}.

\section{Solving the Sudoku puzzle}

\begin{figure}
   \centering
    \includegraphics[width = .9\textwidth]{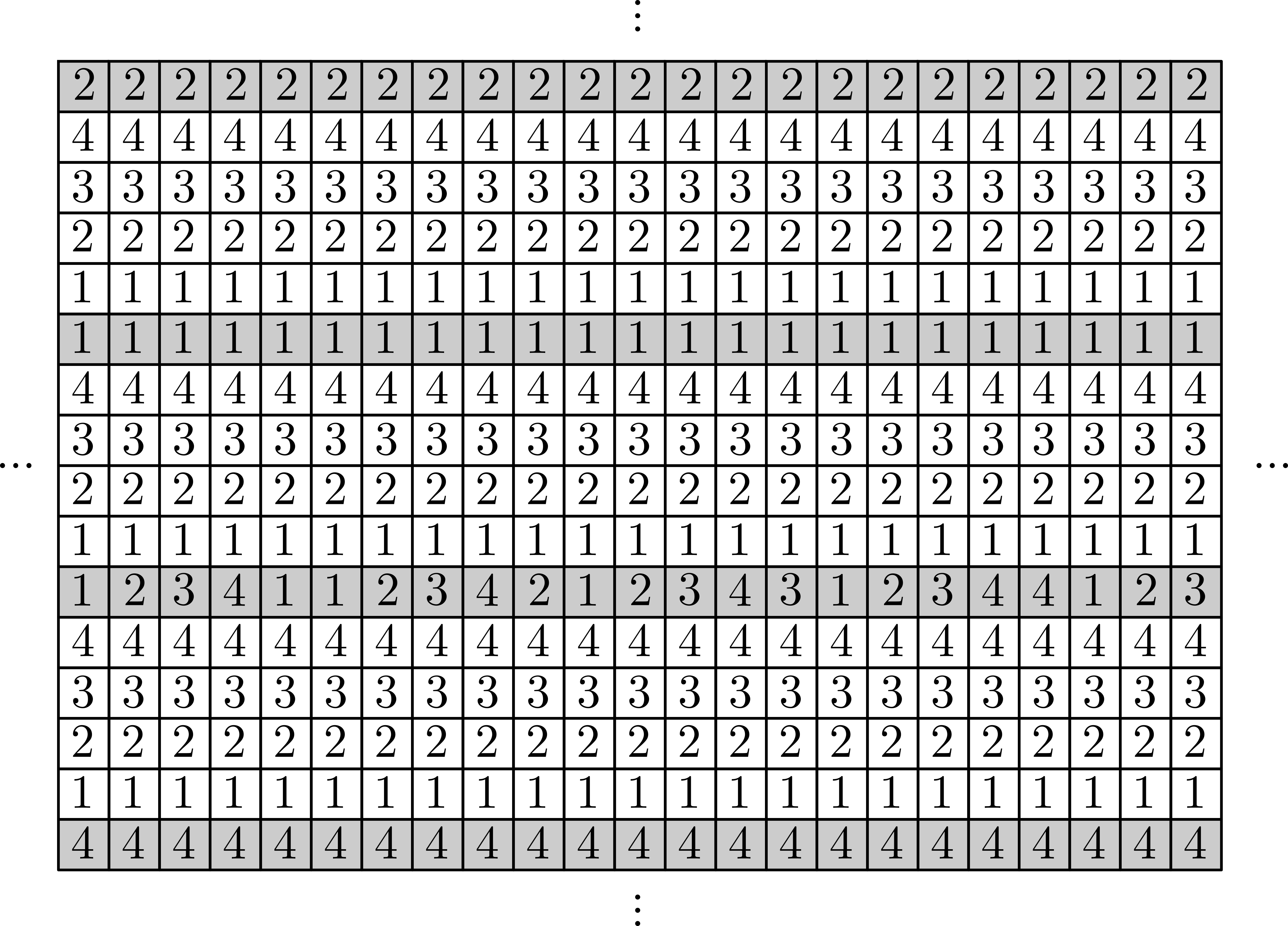}
    \caption{The outcome of applying a normalization to the $\S^2_5$-Sudoku solution in   \figref{fig:sudoku}. Here, all the grey rows are indexed by a multiple of $5$.}
  \label{fig:rescale}
\end{figure}

In this section we briefly sketch how to establish  \thmref{sudoku-aperiodic}.  Observe that the functions in $\S^2_p(\{1,\dots,N\})$ are affine outside of at most one ``bad coset'' of $p\Z$.  Thus, a $\S^2_p$-Sudoku solution $F \colon \{1,\dots,N\} \times \Z \to (\Z/p\Z)^\times$ is ``almost affine'' on every non-horizontal line $\{ (n, jn+i): n=1,\dots,N\}$, in the sense that there are constants $A_{i,j}, B_{i,j} \in \Z/p\Z$, not both zero, such that
$$ F(n,jn+i) = A_{i,j} n + B_{i,j}$$
whenever $n \in \{1,\dots,N\}$ is such that $ A_{i,j} n + B_{i,j} \neq 0$.  It turns out, after some combinatorial case-checking, that (for $p$ large enough) this is enough to force $F$ to behave in an almost affine manner on all of $\{1,\dots,N\} \times \Z$. More precisely, one can establish that there exist coefficients $A,B,C \in \Z/p\Z$, not all zero, such that
\begin{equation}\label{paffine}
F(n,m) = An + Bm + C
\end{equation}
whenever $(n,m) \in \{1,\dots,N\} \times \Z$ is such that $An+Bm+C \neq 0$.  After applying some affine changes of variable, one can normalize to the situation $A=C=0$, $B=1$, so that we now have
$$ F(n,m) = m$$
whenever $m \neq 0 \mod p$; see   \figref{fig:rescale}.  This completely determines the $\S^2_p$-Sudoku solution $F$ as having constant rows outside the coset of ``bad rows'' $\{1,\dots,N\} \times p\Z$. In particular, each of the columns $m\mapsto F(n,m)$ exhibits  $p$-adic structure outside of one coset of $p\Z$.

\begin{figure}
   \centering
    \includegraphics[width = .9\textwidth]{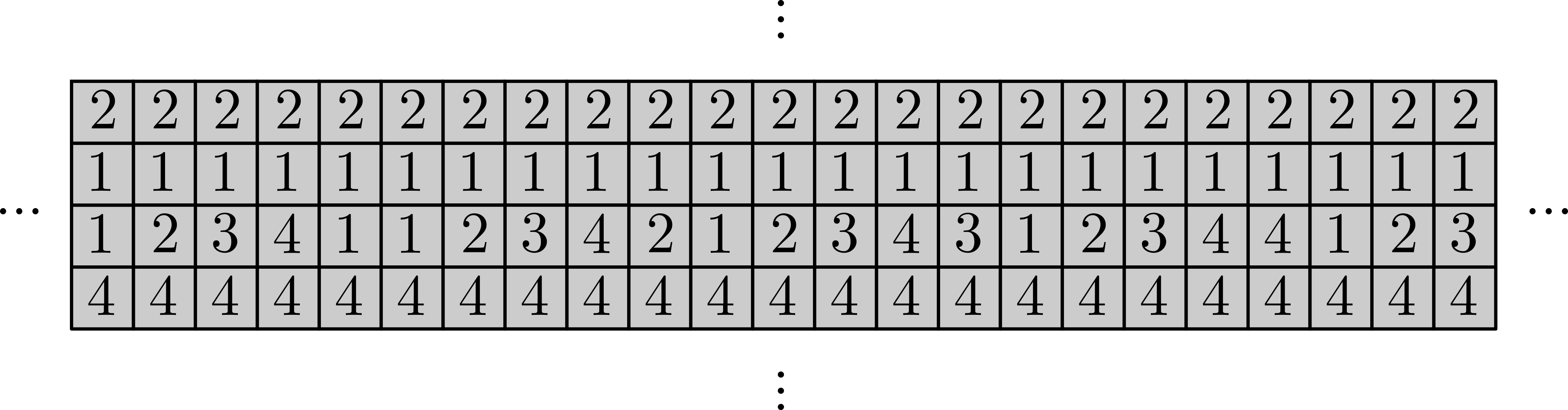}
    \caption{The outcome of applying a tetris move to the $\S^2_5$-Sudoku solution in   \figref{fig:rescale}.}
  \label{fig:tetris}
\end{figure}

The next step is to perform a ``tetris'' move, removing all the constant rows and focusing on the remaining ``bad'' rows, and more precisely by studying the function $F_1 \colon \{1,\dots,N\} \times \Z \to (\Z/p\Z)^\times$ defined by
\begin{equation}\label{tetris}
F_1(n,m) \coloneqq  F(n, pm);
\end{equation}
see   \figref{fig:tetris}.

It turns out that $F_1$ is also a $\S^2_p$-Sudoku solution, and one can repeat the preceding analysis to also conclude that there is also an almost affine relationship
$$ F_1(n,m) = A_1n + B_1m + C_1$$
for some coefficients $A_1,B_1,C_1 \in \Z/p\Z$ (not all zero), and all $(n,m) \in \{1,\dots,N\} \times \Z$ such that $A_1 n + B_1 m + C \neq 0$.  By comparing this with \eqref{paffine} and  \defref{s2p} one can verify that the $B_1$ coefficient for $F_1$ must match the $B=1$ coefficient for $F$.  To get this crucial matching $B_1=B$ that we have the key equation \eqref{key} outside of a coset of $p^2 \Z$, and not merely outside of a coset of $p\Z$ (for this reason, we introduce  a $\S^2_p$-Sudoku rather than merely a $\S^1_p$-Sudoku).  This matching was achieved under the normalization $A=C=0, B=1$, but it is an easy matter to undo the normalization and obtain an analogous matching conclusion $B_1=B$ in general. In particular, we obtain that each of the columns $m\mapsto F(n,m)$ exhibits  $p$-adic affine structure outside of one coset of $p^2\Z$. It is then possible to iterate this matching and conclude that each of the columns $m \mapsto F(n,m)$ exhibits $p$-adic structure without any exception (which basically means that behaves like a rescaled version of $f_p$); as the columns are non-constant, they must therefore be non-periodic, giving the claim.



\begin{thebibliography}{99999}
    		
    			\bibitem[B66]{Ber}
		 R. Berger, \emph{The undecidability of the domino problem}, Memoirs of the American Mathematical Society,    \textbf{66}  (1966) 72.
		 
		\bibitem[B64]{Ber-thesis}
		R. Berger. The Undecidability of the Domino Problem. PhD thesis, Harvard University, 1964.

\bibitem[B26]{besicovitch}
A.S. Besicovitch, \emph{On generalized almost periodic functions}, Proc. London Math. Soc., \textbf{25} (1926), 495--512.
		
		\bibitem[B20]{BH}
		S. Bhattacharya,  \emph{Periodicity and Decidability of Tilings of $\Z^2$}, 
		Amer. J. Math., \textbf{142}, (2020),  255--266.
		
		\bibitem[BN91]{bn}
		D. Beauquier, M. Nivat, \emph{On translating one polyomino to tile the plane}, Discrete Comput. Geom. \textbf{6} (1991), no. 6, 575--592.


         
         \bibitem[C96]{culik}
		K. Culik II., 
		\emph{An aperiodic set of 13 Wang tiles},
		 Discrete Math., \textbf{160} (1996),  245--251.
		 
		 \bibitem[G-BN91]{gbn}
D. Girault-Beauquier, M. Nivat, \emph{Tiling the plane with one tile}, Topology and category theory in computer science (Oxford, 1989), 291–333, Oxford Sci. Publ., Oxford Univ. Press, New York, 1991.
         
         		
         			 \bibitem[GT20]{GT}
		R. Greenfeld, T. Tao, 
		\emph{The structure of translational tilings in $\Z^d$}, Discrete Analysis (2021):16, 28 pp.
	
    		
    			 \bibitem[GT21]{GT2}
    		R. Greenfeld, T. Tao, 
    		\emph{Undecidable translational tilings with only two tiles, or one nonabelian tile}, {\tt arXiv:2108.07902}.

    			 \bibitem[GT22]{greenfeld-tao-ptc}
    		R. Greenfeld, T. Tao, 
    		\emph{A counterexample to the periodic tiling conjecture}, in preparation.
      
      \bibitem[K96]{kari}
J. Kari,  \emph{A small aperiodic set of Wang tiles}, Discrete Math., \textbf{160} (1996), 259--264.

\bibitem[Ken92]{ken}
R. Kenyon, \emph{Rigidity of planar tilings}, Invent. Math., \textbf{107} (1992), 637--651.
    		
    		\bibitem[L95]{L}
    		J. C. Lagarias, \emph{Meyer's Concept of Quasicrystal
and Quasiregular Sets}, Commun Math Phys \textbf{179}, 365--376 (1996)
    		
    		
    		
    		 \bibitem[LW96]{LW}
    J. C. Lagarias, Y. Wang, \emph{Tiling the line with translates of one tile}, Invent. Math. \textbf{124} (1996), no. 1-3, 341--365.
    
      \bibitem[M70]{M70}
    Y. Meyer, \emph{Nombres de Pisot, nombres de Salem et analyse harmonique}, Lecture Notes in
Mathematics \textbf{117} (1970), Springer-Verlag.
    
    \bibitem[M95]{M95}
    Y. Meyer, \emph{Quasicrystals, diophantine approximation and algebraic numbers}, Axel, F., Gratias, D. (eds) Beyond Quasicrystals. Centre de Physique des Houches, vol 3. Springer, Berlin, Heidelberg (1995).
    
    \bibitem[MSS22]{bgu}
    T. Meyerovitch, S.  Sanadhya, Y. Solomon,
    \emph{A note on Reduction of tiling problems}, in preparation.
    
  
    
    			\bibitem[N77]{N}
		D. J. Newman, \emph{Tesselation of integers}, J. Number Theory \textbf{9} (1977), no. 1, 107--111.
		
		\bibitem[O09]{ollinger}
N. Ollinger, \emph{Tiling the plane with a fixed number of polyominoes},
Lecture Notes in Comput. Sci., 5457, Springer, Berlin, 2009.

	\bibitem[S98]{szegedy}
	M. Szegedy, \emph{Algorithms to tile the infinite grid with finite clusters}, Proceedings of the
39th Annual Symposium on Foundations of Computer Science (FOCS ’98), IEEE Computer Society, Los Alamitos, CA, (1998), 137--145.

\bibitem[S06]{S}
B. Solomyak, \emph{Lectures on Tilings and Dynamics}, EMS Summer School on Combinatorics, Automata, and Number Theory, Liege, May 2006, \href{https://u.math.biu.ac.il/~solomyb/RESEARCH/notes6.pdf}{link}.
    		
    		\bibitem[W75]{wang}
		 H. Wang, \emph{Notes on a class of tiling problems}, Fundamenta Mathematicae, \textbf{82} (1975), 295--305.
    		
    	\end{thebibliography}
    \end{document}